\documentclass[a4paper,12pt]{article}
\usepackage{bm,amsmath,amsfonts,amssymb}%,isolatin1 % for accents etc. not used
\date{24 April 2002}%6 Dec 2001}
\title{Generalized function algebras as sequence space algebras}
\makeatletter
\author{Antoine Delcroix\thanks{E-mail: \tt Antoine.Delcroix@univ-ag.fr}
\and Maximilian F. Hasler\thanks{E-mail: \tt
Maximilian.Hasler@martinique.univ-ag.fr}
\and Stevan Pilipovi\'c\thanks{E-mail: \tt pilipovic@im.ns.ac.yu}
\and Vincent Valmorin\thanks{E-mail: \tt Vincent.Valmorin@univ-ag.fr}}

\catcode`\·\active \def·{\cdot}

\def\CC{{\cal C}}
\def\C{{\mathbb C}}

\def\CEP{\ensuremath{(\mathcal C,\mathcal E, \mathcal P)}}

\def\F{{\cal F}}
\def\G{{\cal G}}

\def\INT#1d{\int#1\,\mathrm d}

\def\ind{\operatorname*{ind\,lim\vphantom p}}
\def\impl{\mathop{~\Longrightarrow~}}
\def\K{{\cal K}}

\def\ola#1{\rlap{$\displaystyle\overleftarrow{\phantom{#1}}$}#1}
\def\olra{\overleftrightarrow}
\def\ora{\overrightarrow}

\def\Mid{~\Big|~}
\def\N{{\mathbb N}}
\def\nb#1{\noindent{\bf #1}}

\def\p#1{\left( #1 \right)}
\def\R{{\mathbb R}}
\def\supp{\operatorname*{supp}}
\def\set#1{\left\{\,#1\,\right\}}
\def\ultra#1{|\!|\!|\,#1\,|\!|\!|} % |||\([^|]*\)|||

\def\veps{\varepsilon}
\def\Z{{\mathbb Z}}

\newtheorem{proposition}{Proposition}[section]
\newtheorem{propdef}[proposition]{Proposition--Definition}
\newtheorem{example}{Example}
\newtheorem{remark}{Remark}
\newenvironment{proof}[1][Proof]{\par\nb{#1.} }
{ \mbox{}\hfill$\square$\par\medskip}

\advance\textheight\topmargin
\begin{document}
\maketitle

\begin{abstract}
A topological description of various generalized function algebras over
corresponding basic locally convex algebras is given.  The framework
consists of algebras of sequences with appropriate
ultra(pseudo)metrics defined by sequences of exponential weights.
Such an algebra with embedded Dirac's delta distribution induces
discrete topology on the basic space. This result is in analogy to
Schwartz' impossibility result concerning multiplication of
distributions.%
\def\thefootnote{}\footnote{~\\[-1\baselineskip]
MSC: % http://www.ams.org/msc/46Axx.html
46A45 (Sequence spaces), % (including Köthe sequence spaces) [See also 46B45]
46F30 (Generalized functions for nonlinear analysis).\\
	% (Rosinger, Colombeau, nonstandard, etc.)
secondary:
46E10, % Topological linear spaces of continuous,
	% differentiable or analytic functions
%46A03, General theory of locally convex spaces
46A13, % Spaces defined by inductive or projective limits (LB, LF, etc.)
46A50, % Compactness in topological linear spaces; angelic spaces, etc.
% 46H05 General theory of topological algebras
46E35, % Sobolev spaces and other spaces of ``smooth'' functions,
% embedding theorems, trace theorems
% FOR SECOND PART :
% 46A22 extension and lifting of functionals and operators [See also 46M10]
46F05.% Topological linear spaces of test functions,
% distributions and ultradistributions [See also 46E10, 46E35]
%46F10 Operations with distributions
%46F15 Hyperfunctions, analytic functionals
%46F20 Distributions and ultradistributions
% as boundary values of analytic functions [See also 30D40, 30E25, 32A40]
}%
\end{abstract}

\section{Introduction}

%%% NEW
After Schwartz' ``impossibility result'' \cite{schw} for algebras of
generalized functions with prescribed list of (natural) assumptions,
several new approaches had appeared with the aim of applications in
non linear problems. We refer to the recent monograph~\cite{GKOS} for
the historical background as well as for the list of relevant references
mainly for algebras of generalized functions today called Colombeau
type algebras (see~\cite{bia,BC,co,gfks,ob}). 
% [1],[2],[3],[7],[9]...). 
Colombeau and all other successors introduce algebras of generalized
functions
% by algebraic methods
%We give a topological description of Colombeau type
%algebras~\cite{bia,BC,co,gfks,ob}, 
%originally introduced by Colombeau
through purely algebraic methods. By now, these algebras have become an
important tool in the theory of PDE's, stochastic analysis, differential
geometry and general relativity.  We show that such algebras fit in the
general theory of well known sequence spaces forming appropriate algebras.
These classes of algebras of sequences are simply determined by a locally
convex algebra $E$, and a sequence of weights (or sequence of sequences),
which serve to construct ultra(pseudo)metrics.

%%% NEW %%% e-mail sometimes changes ``From'' to ``>From''
{}From the beginning, the topological questions concerning such algebras
were important. We refer to the papers~\cite{BC} where the classical
topology
% (I forgot the name given by scarpalezos for this topology)
and a uniform structure were introduced in order to consider
generalized functions as smooth functions in appropriate quotient spaces.
Then was introduced the sharp topology~\cite{O}
% (nowdays called the sharp one) 
in connection with the well posedness of the Carleman system with
measures as initial data. Later it was independently
reintroduced and analyzed in~\cite{Scarp}, where
the name ``sharp topology'' appeared.
%But there was always the question
It remained an open question whether the introduced topologies were
``good enough'', because they induced always the discrete topology on
the underlying space.

We show that the topology of a Colombeau type algebra containing
Dirac's delta distribution $\delta$ as an embedded Colombeau
generalized function, must induce the discrete topology on the basic
space $E$. This is an analogous result to Schwartz' ``impossibility
result'' concerning the product of distributions (cf.  Remark at the
end of the Note).

We mention that distribution,
ultradistribution and hyperfunction type spaces can be embedded in
corresponding algebras of sequences with exponential weights
(cf.~\cite{dhpv2}).
More general concepts of generalized functions not anticipating
embeddings as well as regularity properties of generalized functions
can be found in~\cite{ros} and~\cite{CEP}. Another interesting approach
is given in~\cite{todor}

To be short, we give most examples only for spaces of functions defined on
$\R^s$, although the generalization to an open subset of $\R^s$ is
straightforward.

\section{General construction}

Consider a positive sequence $ r={(r_n)}_n\in(\R_+)^\N $ decreasing to
zero.  If $p$ is a seminorm on a vector space $E$, we define for
$f={(f_n)}_n\in E^\N$
$$
    \ultra f_{p,r} = \limsup_{n\to\infty} \big( p(f_n) \big)^{r_n}
$$
with values in $\overline\R_+ = \R_+\cup\set\infty$.
Denote $\tilde E^\N=\{f\in E^\N \mid \ultra f_{p,r} <\infty\}.$
\\
Let $\left( E_\nu^\mu,p_{\nu\,}^\mu\right)_{\mu,\nu\in\N}$ be
a family of semi-normed algebras over $\R$ or $\C$ such that
\begin{align*}
	\forall\mu,\nu\in\N :~ E_\nu^{\mu+1} \hookrightarrow E_\nu^\mu
\text{ ~ and ~ }
	E_{\nu+1}^\mu \hookrightarrow E_\nu^\mu
\text{~ (resp. }
	E_\nu^\mu \hookrightarrow E_{\nu+1}^\mu ~)~,
\end{align*}
where $\hookrightarrow$ means continuously embedded.
~(For the $\nu$ index we consider inclusions in the two directions.)~
%
%This implies that there exist constants $ C_\nu^\mu\in\R _+ $
%such that
%$ \forall\mu,\nu\in\N : p_\nu^\mu\leq C_\nu^\mu\,p_{\nu+1}^\mu $
%(on domains where both seminorms are defined),
%but without loss of generality one can take these constants equal to 1.
%The same holds for the $_\nu$ indices.
%
Then let
$\displaystyle 
        \ola{E}
        =\projlim_{\mu\to\infty} \projlim_{\nu\to\infty} E_\nu^\mu
        =\projlim_{\nu\to\infty} E_\nu^\nu
$, (resp.
$\displaystyle
        \ora{E}
        =\projlim_{\mu\to\infty} \ind_{\nu\to\infty}E_\nu^\mu)
$.
Such projective and inductive limits are usually considered with norms
instead of seminorms, and with the additional assumption that in the
projective case sequences are reduced, while in the inductive case for
every $\mu\in\N$ the inductive limit is regular, i.e. a set
$A\subset{\ind\limits_{\nu\to\infty}}E_\nu^\mu$ is bounded iff
it is contained in some $E^\mu_\nu$ and bounded there.

Define (with $p\equiv\p{p_\nu^\mu}_{\nu,\mu}$)
\begin{align*}
	\ola\F_{p,r} &= \set{ f\in\ola{E}^\N \Mid
	\forall\mu,\nu\in\N : \ultra f_{p_\nu^\mu,\,r} < \infty } ~,
\\
	\ola\K_{p,r} &= \set{ f\in\ola{E}^\N \Mid
	\forall\mu,\nu\in\N : \ultra f_{p_\nu^\mu,\,r} = 0 }
\\
\text{(resp.}~~
	\ora\F_{p,r} &= \bigcap_{\mu\in\N} \ora\F_{p,r}^\mu
~,~~	\ora\F_{p,r}^\mu = \bigcup_{\nu\in\N}
	\set{ f\in\p{ E_\nu^\mu }^\N \Mid \ultra f_{p_\nu^\mu,r} < \infty }
~,\\
	\ora\K_{p,r} &= \bigcap_{\mu\in\N} \ora\K_{p,r}^\mu
~,~~	\ora\K_{p,r}^\mu = \bigcup_{\nu\in\N}
	\set{ f\in\p{ E_\nu^\mu }^\N \Mid \ultra f_{p_\nu^\mu,r} = 0 }
~\text)~.
\end{align*}
%
%Then, we have the following
%
\begin{propdef}~\parskip-1ex
\begin{enumerate}
\item[(i)] Writing $\olra\cdot$ for both, $\ola\cdot$ or
$\ora\cdot$, we have that $\olra\F_{p,r}$ is an algebra and
$\olra\K_{p,r}$ is an ideal of $\olra\F_{p,r}$; thus, $\olra\G_{p,r}=
\olra\F_{p,r}/\olra\K_{p,r}$ is an algebra.

\item[(ii)] For every $\mu,\nu\in\N,~ d_{p_\nu^\mu}: {(E_\nu^\mu)}^\N
\times {(E_\nu^\mu)}^\N \to\overline\R_+$ defined by
$d_{p_\nu^\mu}(f,g)= \ultra{f-g}_{p_\nu^\mu,r}$ is an
ultrapseudometric on ${(E_\nu^\mu)}^\N$.  Moreover,
$(d_{p_\nu^\mu})_{\mu,\nu}$ induces a topological algebra\footnote[1]
{\label{a} over $(\C^\N,\ultra\cdot_{|\cdot|})$, not over $\C$:
scalar multiplication is not continuous because of (\protect\ref3)}
structure on $\ola\F_{p,r}$
% (since $d_{p_\nu^\mu}(0,f\cdot g)\leq
% d_{p_\nu^\mu}(0,f)\,d_{p_\nu^\mu}(0,g)$) 
such that the intersection of
the neighborhoods of zero equals $\ola\K_{p,r}$.

\item[(iii)] From (ii), $\ola\G_{p,r}= \ola\F_{p,r}/\ola\K_{p,r}$
becomes a topological algebra (over generalized numbers
$\C_r=\G_{|\cdot|,r}$) whose topology can be defined by the family of
ultrametrics $(\tilde{d}_{p_\nu^\mu})_{\mu,\nu}$ where
$\tilde{d}_{p_\nu^\mu} ([f],[g])=d_{p_\nu^\mu}(f,g)$, $[f]$ standing
for the class of $f$.

\item[(iv)] If $\tau_\mu$ denote the inductive limit topology on
$\F_{p,r}^\mu=\bigcup_{\nu\in\N}((\tilde{E_\nu^\mu})^\N,d_{\mu,\nu})$,
$\mu\in\N$, then $\ora\F_{p,r}$ is a topological algebra\footnotemark[1%\ref{a}
]
for the
projective limit topology of the family $(\F_{p,r}^\mu,\tau_\mu)_\mu$.

\end{enumerate}
\end{propdef}

\begin{proof}We use the following properties of $\ultra\cdot$:
\begin{align}
\label{1} \forall x,y\in E^\N: \ultra{ x+y } &\le \max(\ultra x,\ultra y) ~,\\
\label{2} \forall x,y\in E^\N, \ultra{ x·y } &\le \ultra x·\ultra y ~,\\
\label{3} \forall \lambda\in\C^*,x\in E^\N: \ultra{\lambda\,x} &= \ultra x ~,
\end{align}
%The second follows directly from $p(x_n·y_n)=p(x_n)·p(y_n)$,
They are consequences of basic properties of seminorms and of
$p(x_n+y_n)^{r_n}\le2^{r_n}\max(p(x_n),p(y_n))^{r_n}$ for (\ref1).
Using the above three inequalities, (i)--(iv) follow straightforwardly from
the respective definitions.
\end{proof}

%% \subsection{Examples}

%\subsubsection{Colombeau generalized numbers}

\begin{example}{\sc(Colombeau generalized numbers and ultracomplex numbers)}
Take $E_\nu^\mu=\R$ or $\C$, and $p_\nu^\mu=|\cdot|$ (absolute value)
for all $\mu,\nu\in\N$.  Then, for $r_n=\frac1{\log n}$, we get the
ring of Colombeau's numbers $\overline\R$ or $\overline\C$.

\noindent With the sequence $r_n=n^{-1/m}$ for some fixed $m>0$,
we obtain rings of ultracomplex numbers $\overline\C^{p!^m}$
(cf.~\cite{dhpv2}).

\end{example}

% \begin{remark}
% To be precise, $\ola\G_{p,r}$ (resp. $\olra\F_{p,r}$) are
% topological algebras over generalized numbers $\C_r=\G_{|\cdot|,r}$
% (resp.  $\F_{|\cdot|,r}$) with $d_{|\cdot|,r}$ topology.  Scalar
% multiplication by $\lambda\in\C$ is not continuous in $\lambda=0$
% because $\forall\lambda\ne0:\ultra{\lambda\,f}=\ultra f$.
% \end{remark}

\begin{example}\label{Ex-Sobolev}
Consider a Sobolev space $E=W^{s,\infty}(\Omega)$
for some $s\in\N $. The corresponding Colombeau type algebra is
defined by $\G_{W^{s,\infty}}=\F/\K$, where
\begin{align*}
   \F &= \set{ u\in (W^{s,\infty})^\N \mid
   \limsup \left\| u_n \right\|_{s,\infty\strut}^{\frac1{\log n}} < \infty } ~,
\\ \K &= \set{ u\in (W^{s,\infty})^\N \mid
   \limsup \left\| u_n \right\|_{s,\infty\strut}^{\frac1{\log n}} = 0 } ~.
\end{align*}
\end{example}

\begin{example}{\sc(simplified Colombeau algebra)}\label{Ex-Col-Simpl}
Take $ E_\nu^\mu = \CC^\infty (\R^s) $,
$$
  p_\nu^\mu(f) = \sup_{|\alpha|\leq \nu,\,|x|\leq\mu} |f^{(\alpha)}(x)| ~,
$$
and $r=\frac{1}{\log}$. Then, $\ola\G_{p,r} = \ola\F_{p,r} /\ola\K_{p,r} $
is the simplified Colombeau algebra.\\
With slight modifications, the full Colombeau algebra can also be described
in this setting.
\end{example}

\section{Completeness}

Without assuming completeness of $ \olra E $, we have\vskip-3ex
\begin{proposition}
\begin{enumerate}
\item[(i)] $\ola\F_{p,r}$ is complete.
\item[(ii)] If for all $\mu\in\N$, a subset of $\ora\F^\mu_{p,r}$ is
bounded iff it is a bounded subset of $\left(E_\nu^\mu\right)^\N$ for
some $\nu\in\N$, then $\ora\F_{p,r}$ is sequentially complete.
\end{enumerate}
\end{proposition}

\begin{proof}
If ${(f^m)}_{m\in\N}\in\ola\F_{p,r}$ is a Cauchy sequence,
there exists a strictly increasing sequence ${(m_\mu)}_{\mu\in\N}$
of integers such that
$$
	\forall\mu\in\N ~~ \forall\,k,\ell\ge m_\mu:
	\limsup_{n\to\infty} p^\mu_\mu\p{f^k_n-f^\ell_n}^{r_n}<\frac1{2^\mu} ~.
$$
Thus, there exists a strictly increasing sequence ${(n_\mu)}_{\mu\in\N}$
of integers such that
$$
	\forall\mu\in\N ~~ \forall\,k,\ell\in[m_\mu,m_{\mu+1}] ~~
	\forall n\ge n_\mu: p^\mu_\mu\p{f^k_n-f^\ell_n}^{r_n}<\frac1{2^\mu} ~.
$$
(Restricting $k,\ell$ to $[m_\mu,m_{\mu+1}]$ allows to take
$n_\mu$ independent of $k,\ell$.)\\
\newcommand\mm{{\bm\mu}}
\newcommand\sm[1]{^{\smash{m_{#1}}}}
Let $ \mm(n)=\sup\set{\mu\mid n_\mu\le n}$, and
consider the diagonalized sequence
$$
	\bar f={(f\sm{\mm(n)}_n)}_n
\text{ , i.e. ~ }
	\bar f_n=\begin{cases}
			f^{m_0}_n & \text{if ~ } n\in[n_0,n_1)\\...\\
			f\sm\mu_n & \text{if ~ } n\in[n_\mu,n_{\mu+1})\\
		...\end{cases} ~.
$$
Now let us show that
%$\lim_{m\to\infty}\ultra{f^m-\bar f}=0$.
$f^m\to\bar f$ in $\ola\F_{p,r}$, as $m\to\infty$.
Indeed, for $\veps$ and $p^{\mu_0}_\nu$ given,
choose $\mu>\mu_0,\nu$ such that $\frac1{2^\mu}<\frac12\veps$.
As  $p^\mu_\nu$ is increasing in both indices, 
we have for $m>m_\mu$ ~(say $m\in[m_{\mu+s},m_{\mu+s+1}]$)~:
\begin{eqnarray*}
	p^{\mu_0}_\nu(f^m_n-\bar f_n)^{r_n}
	&\le& p^\mu_\mu(f^m_n-f_n\sm{\mm(n)})^{r_n}
\\	&\le& p^\mu_\mu(f^m_n-f_n\sm{\mu+s+1})^{r_n}
	+ \sum_{\mu'=\mu+s+1}^{\mm(n)-1}
		p^{\mu'}_{\mu'}(f\sm{\mu'}_n-f\sm{\mu'+1}_n)^{r_n}
\end{eqnarray*}
and for $n>n_{\mu+s}$, we have
% $\mm(n)>\mu+s$
of course $n\ge n_{\mm(n)}$, thus finally
\begin{eqnarray*}
	p^{\mu_0}_\nu(f^m_n-\bar f_n)^{r_n}
	< \sum_{\mu'=\mu+s}^{\mm(n)} \frac1{2^{\mu'}}
	< \frac2{2^\mu} < \veps
\end{eqnarray*}
and therefore $ f^m\to\bar f$ in $ \ola\F $.

%% END NEW ADDITIONS

For $\p{f^m}_m$ a Cauchy net in $\ora\F_{p,r}$, the proof
requires some additional considerations.
We know that for every $\mu$ there is $ \nu(\mu) $ such that
$$
  p_{\nu(\mu)}^\mu \left( f_n^m - f_n^p \right) ^{r_n}
  <\veps_\mu,
$$
where ${(\veps_\mu)}_\mu$ decreases to zero.
For every $\mu$ we can choose $\nu(\mu)$ so that
$ p_{\nu(\mu)}^\mu  \le  p_{\nu(\mu+1)}^{\mu+1}$.
Now by the same arguments as above, we prove the completeness
in the case of $\ora\F_{p,r}$.
\end{proof}

%%%%%%%%%%%%% EMBEDDINGS %%%%%%%%%%%%%%%%

\section{General remarks on embeddings of duals}\label{embedding-remark}

Under mild assumptions on $\olra E$, we show that our algebras of
(classes of) sequences contain embedded elements of strong dual spaces
$\olra E'$. First we consider the embedding of the delta distribution.
We show that general assumptions on test spaces or on a delta
sequence lead to the non-boundedness of a delta sequence in $\olra E$.

%Denote by $\CC^0(\R^s)$ the space of continuous functions with
%projective topology given by sup norms on the balls of radius
%$\nu\in\N^*$: $ p_\nu(f) = \sup \set{ |f(x)| ; |x|\le\nu } $.

We consider $F=\CC^0(\R^s)$, the space of continuous functions with the
projective topology given by sup norms on the balls $B(0,n),~n\in\N^*$,
or $F=\K(\R^s)=\ind_{n\to\infty}(\K_n,\|\cdot\|_\infty)$, where
$$ \K_n = \set{ \phi\in\CC(\R^s) \mid \supp\psi\subset B(0,n) }~.$$
(Recall, $\K'(\R^s)$ is the space of Radon measures.)

%We assume in the sequel that $\olra{E}$ is a dense subspace of
%$\CC^0(\R^s)$ and the inclusion mapping $\olra{E}\to\CC^0(\R^s)$ is
%continuous.

We assume that $\olra E$ is dense in $F$ and $\olra E\hookrightarrow F$. 
This implies that $\delta\in F'\subset \olra E'$.

\begin{proposition}\label{prop10}~\begin{enumerate}
\item[(i)] For $F=\CC^0(\R^s)$, a sequence ${(\delta_n)}_n$ with elements in
$\olra E\cap(\CC^0(\R^s))'$ such that
$$
	\exists M>0, \forall n\in\N: \sup_{|x|>M} |\delta_n(x)|<M ~,
$$
converging weakly to $\delta$ in $\olra E'$, cannot be bounded in $\olra E$.

With $F=\K(\R^s)$, if for $\p{\delta_n}_n\in(\olra E)^\N$
there exists a compact set $K$ so that $\supp\delta_n\subset K,~n\in\N^*$,
then this sequence can not be bounded in $\olra E$.

\item[(ii)] Assume:\\
1. Any $\phi\in\olra E$ defines an element of $F'$  by
$\psi\mapsto\INT_{\R^s}\phi(x)\,\psi(x)dx$.\\
2. If ${(\phi_n)}_n$ is a bounded sequence in $\olra E$, then 
$\sup\limits_{n\in\N,x\in\R^s}|\phi_n(x)|<\infty$.

Then, if $\olra{E}$ is sequentially weakly dense in $\olra{E}'$ and
${(\delta_n)}_n$ is a sequence in $\olra{E}$ converging weakly to
$\delta$ in $\olra E'$, then ${(\delta_n)}_n$ can not be bounded in
$\olra E$.
\end{enumerate}
\end{proposition}
% Let us remark that the technical condition on $\p{\delta_n}_n$ is
% satisfied for the usual choices of $\delta$--nets.

\begin{proof}
(i)
We will prove the assertion only for $F=\CC^0(\R^s).$
Let us show that ${(\delta_n)}_n$ is not bounded in $\olra E$.
First consider $\ola{E}$.
Boundedness of $\p{\delta_n}_n$ in $\ola{E}$ implies:
$\forall\mu\in\N,~\forall\nu\in\N,~\exists C_1>0,~\forall n\in\N:
p_\nu^\mu(\delta_n)<C_1$.
Continuity of $\ola{E}\hookrightarrow\CC^0(\R^s)$ gives
\[
	\forall k\in\N, \exists\mu\in\N, \exists\nu\in\N, \exists C_2>0,
	\forall\psi\in\ola E :
	\sup_{|x|<k} |\psi(x)| \leq C_2\,p_\nu^\mu(\psi) ~.
\]
It follows:
$\exists C>0,\forall n\in\N:\sup_{x\in\R^s} |\delta_n(x)|<C$,
which is impossible. To show this, take $\psi\in\CC^0(\R^s)$
so that it is positive, $\psi(0)=C+1$, and $\int\psi<1$.
The assumption $\delta_n\in(\CC^0(\R^s))'$ implies that it
acts on $\CC^0(\R^s)$ by $\psi\mapsto\INT\delta_n(x)\,\psi(x)dx$.
This gives
$
	C+1=\psi(0)\gets\left|\INT\delta_n\psi dx\right| \le C
$.\\
For $\ora{E}$, simply exchange $\forall\nu\leftrightarrow\exists\nu$
in the above.
 
(ii) Assumption 2. and the boundedness of $\p{\delta_n}_n$ in
$\olra{E}$ would imply: $\exists C>0,\forall n\in\N:\sup_{x\in\R^s}
|\delta_n(x)|<C$. Now by assumption 1. we conclude the proof as in
part (i).
\end{proof}

\begin{remark}
One can take for $\olra E$ one of Schwartz' test function spaces or
Beurling or Roumieau test function space of ultradifferentiable
functions.  Since the delta distribution lives on all functions which
are continuous at zero, one can consider also $F$ and $\olra E$ to
consist of holomorphic functions with appropriate topologies. This was
the reason for considering $\CC^0$, although there are many classes of
test spaces which would imply the necessary accommodation of conditions
of the previous assertion.
\end{remark}

Thus, the appropriate choice of a sequence $r$ decreasing to 0 appears
to be important to have at least $\delta$ embedded into the
corresponding algebra. It can be chosen such that for all $\mu\in\N$
and all $\nu\in\N$ (resp. some $\nu\in\N$ in $\ora{E}$ case),
$\limsup_{n\to\infty}p_\nu^\mu(\delta_n)^{r_n} = A_\nu^\mu$ and
$\exists\mu_0,\nu_0: A_{\nu_0}^{\mu_0} \ne 0$.

%\begin{itemize}
%\item  in $\ola{E}$ case, $\forall\mu,\nu\in\N $,
%$ %\begin{equation}
%  \limsup_{n\to\infty}p_\nu^\mu\p{ \delta_n }^{r_n} =
%  A_\nu^\mu
%\text{~ and ~}
%  \exists \mu_0, \nu_0: A_{\nu_0}^{\mu_0} \neq 0
%%% \tag{1}\label{(1)}
%$.%\end{equation}

%%% and $\exists\mu_{0},\nu_{0}$ such that $A_{\nu_{0}}^{\mu_{0}}\neq0$.

% \item in $\ora{E}$ case, $\forall\mu\in\N $,
% $\exists\nu\in\N $ such that the previous equation holds.
% \end{itemize}

So the  embedding of duals into corresponding algebras is realized
on the basis of two demands:
\begin{enumerate}
\item $\olra E$ is weakly sequentially dense in $\olra E'$.

\item There exists a sequence $\p{r_n}_n$ decreasing to zero, such
that for all $f\in\olra{E}'$ and corresponding sequence $\p{f_n}_n$ in
$\olra E$, $f_n\to f$ weakly in $\olra E'$, we have for all $\mu$ and
all $\nu$ (resp. some $\nu$), $\limsup\limits_{n\to\infty}
p_\nu^\mu(f_n)^{r_n}<\infty$.
\end{enumerate}

\begin{remark}
In the definition of sequence spaces $\olra\F_{p,r}$, we assumed
$r_n\searrow0$ as $n\to\infty$.
In principle, one could consider more general sequences of weights.
For example, if $r_n\in(\alpha,\beta)$, $0<\alpha<\beta$, then $\olra
E$ can be embedded, in the set-theoretical sense, via the canonical
map $f\mapsto(f)_n$ ($f_n=f$).
%
%For beta\le1 the original topology is induced
%
If $r_n\to\infty$, $\olra E$ is no more included in $\olra\F_{p,r}$.

In the case we consider ($r_n\to0$), the induced topology on
$\olra{E}$ is obviously a discrete topology.
But this is necessarily so, since we want to have ``divergent''
sequences in $\olra\F_{p,r}$. Thus, in our construction, in order to have an
appropriate topological algebra containing ``$\delta$'', it is
unavoidable that our generalized topological algebra induces a
discrete topology on the original algebra $\olra{E}$.\\
In some sense, in our construction this is the price to pay, in
analogy to Schwartz' impossibility statement for multiplication of
distributions~\cite{schw}.
\end{remark}

% This Note is the first part of the summary of a bigger paper
% which will be published elsewhere.
% In a second Note, we will deal with questions of embedding spaces
% of distributions, ultradistributions and periodic hyperfunctions
% into the spaces considered here.
%
% Different notions of association will also be discussed there.

\section{Sequences of scales}

The analysis of previous sections can be extended to the case where
we consider a sequence $\p{r^m}_m$ of decreasing null sequences
${(r_n^m)}_n$, satisfying one of the following additional conditions:
\begin{align*}
        \forall m,n\in\N : r_n^{m+1} \ge r_n^{m}  %  \label{I}
\text{ ~ or ~ }
        \forall m,n\in\N : r_n^{m+1} \le r_n^{m} ~. % \label{II}
\end{align*}

Then let, in the first (resp. second) case :
\begin{align*}
  \olra\F_{p,r} &= \bigcap_{m\in\N} \olra\F_{p,r^m} ~,~~
  \olra\K_{p,r} = \bigcup_{m\in\N} \olra\K_{p,r^m} ~,\\
\text{(resp. ~}
  \olra\F_{p,r} &= \bigcup_{m\in\N} \olra\F_{p,r^m} ~,~~
  \olra\K_{p,r} = \bigcap_{m\in\N} \olra\K_{p,r^m}
\text{~ ), ~ where ~}
        p=\left( p_\nu^\mu \right)_{\nu,\mu} ~.
\end{align*}
\begin{proposition}
With above notations, $\olra\G_{p,r}=\olra\F_{p,r}/\olra\K_{p,r}$ is
an algebra.
\end{proposition}
\begin{proof}
In the first case,
$r^{m+1}\ge r^m \impl \ultra f_{r^{m+1}}\ge\ultra f_{r^m}$ if $p(f_n)\ge1$,
hence $\F_{m+1}\subset\F_m$, and conversely, $p(k_n)\le1$ gives $
\ultra k_{r^{m+1}}\le\ultra k_{r^m}$, thus $\K_{m+1}\supset\K_m$.
Thus, 
%intersection for $\F$ and union for $\K$ makes sense, and 
$\F$ is obviously a subalgebra. To see that $\K$ is an ideal, take
$(k,f)\in\K\times\F$. Then $\exists m:k\in\K_m$, but also $f\in\F_m$,
in which $\K_m$ is an ideal. Thus $k\cdot f\in\K_m\subset\K$.\\[1ex]
If $r^m$ is decreasing, $\F_{m+1}\supset\F_m$ and $\K_{m+1}\subset\K_m$.
%, justifying definitions of $\F$ and $\K$.
Because of this inclusion property, $\F$ is a subalgebra.
To prove that $\K$ is an ideal, take $(k,f)\in\K\times\F$, i.e.
$\forall m:k\in\K_m$, and $\exists m':f\in\F_{m'}$.
We need that $\forall m: k·f\in\K_m$.
Indeed,
% So let $m$ be given.\\
if $m\le m'$, then $\K_{m'}\subset\K_m$, thus
$k·f\in\K_{m'}·\F_{m'}\subset\K_{m'}\subset\K_m$,
and if $m'\le m$, then $\F_{m'}\subset\F_m$, thus
$k·f\in\K_m·\F_{m'}\subset\K_m·\F_m\subset\K_m$.
\end{proof}

\begin{example}
$r_n^m=1/|\log a_m(n)| $, where ${(a_m:\N\to\R_+)}_{m\in\Z}$ is
an asymptotic scale, i.e. $\forall m\in\Z:a_{m+1}=o(a_m)$,
$a_{-m}=1/a_m$, $\exists M\in\Z: a_M=o(a_m^2)$. This gives back the
asymptotic algebras of~\cite{ds}.
\end{example}

\begin{example}
Colombeau type ultradistribution and periodic hyperfunction algebras will
be considered in~\cite{dhpv2}.
\end{example}

\begin{example}
$ r^m = \raisebox{0.3ex}{$\chi$}_{[0,m]}$, i.e. $r_n^m=1$ if $n\le m$
and 0 else, gives the Egorov-type algebras~\cite{egor}, where the
``subalgebra'' contains everything, and the ideal contains only
stationary null sequences (with the convention $0^0=0$).
\end{example}

In the case of sequences of scales our second demand
of previous section should read:
``There exist a sequence $\p{\p{r_n^m}_n}_m$ of sequences
decreasing to 0, such that for all $f\in\olra{E}'$ and corresponding
sequence $\p{f_n}_n$ in $\olra E$, $f_n\to f$ weakly in $\olra E'$,
there exists an $m_0$ such that for all $\mu$ and all $\nu$
(resp. some $\nu$),
$\limsup\limits_{n\to\infty} p_\nu^\mu{(f_n)}^{r_n^{m_0}}<\infty$.''

Topological properties for such algebras are a little more complex,
but the ideas of constructing families of ultra(pseudo)metrics are now
clear.  An important feature of our general concept is to show how
various classes of ultradistribution and hyperfunction type spaces can
be embedded in a natural way into sequence space algebras as considered
in Section 2~\cite{dhpv2}.

\bigskip\parindent0pt\parskip1ex\flushleft

Antoine Delcroix : IUFM des Antilles et de la Guyane,\\
Morne Ferret, BP 399, 97159 Pointe \`a Pitre cedex (Guadeloupe, F.W.I.)\\
Tel.: 0590 21 36 21, Fax : 0590 82 51 11,
E-mail: \texttt{Antoine.Delcroix@univ-ag.fr}

Maximilian F. Hasler : Universit\'e des Antilles et de la Guyane,\\
D\'ept Scientifique Interfacultaire, BP 7209,
97275 Schoelcher cedex (Martinique, F.W.I.)\\
Tel.: 0596 72 73 55, Fax : 0596 72 73 62,
E-mail: \texttt{Maximilian.Hasler@martinique.univ-ag.fr}

Stevan Pilipovi\'c : University of Novi Sad, Inst. of Mathematics,\\
Trg D. Obradovi\'ca 4, 21000 Novi Sad (Yougoslavia).\\
Tel.: 00381 21 58 136, Fax : 00381 21 350 458,
E-mail: \texttt{pilipovic@im.ns.ac.yu}

Vincent Valmorin : Universit\'e des Antilles et de la Guyane,
D\'ept de Math\'ematiques,
Campus de Fouillole,  97159 Pointe \`a Pitre cedex (Guadeloupe, F.W.I.)\\
Tel.: 0590 93 86 96, Fax : 0590 93 86 98,
E-mail: \texttt{Vincent.Valmorin@univ-ag.fr}

%{\begin{itemize}\labelsep=2mm\leftskip=-5mm\def\up#1{$^#1$}
%\item[\up a]\end{itemize}}
%\bigskip
%Proofs should be sent to M. Hasler (address see above).

\end{document}